\documentclass[leqno]{amsart}
\usepackage{amsmath}
\usepackage{mathtools}
\usepackage{amssymb}
\usepackage{amsthm}
\usepackage{graphicx}
\usepackage{enumerate}
\usepackage[mathscr]{eucal}
\usepackage{tikz}
\usetikzlibrary{decorations.text,calc,arrows.meta}
\theoremstyle{plain}

\numberwithin{equation}{section}
\usepackage[english]{babel}

\DeclareMathOperator{\Int}{Int}

\newtheorem{theorem}{Theorem}[section]
\newtheorem{cor}[theorem]{Corollary}

\newtheorem{rem}[theorem]{Remark}

\theoremstyle{definition}
\newtheorem{definition}[theorem]{Definition}

\newtheorem{exm}{Example}[section]

\newcommand{\ol}{\overline}
\newcommand{\be}{\begin{equation}}
\newcommand{\ee}{\end{equation}}
\newcommand{\beas}{\begin{eqnarray*}}
\newcommand{\eeas}{\end{eqnarray*}}
\newcommand{\bea}{\begin{eqnarray}}
\newcommand{\eea}{\end{eqnarray}}

\numberwithin{equation}{section}

\usepackage{hyperref}
\hypersetup{
    colorlinks=true,
    linkcolor=blue,
    filecolor=magenta,      
    urlcolor=cyan,
}
\usepackage[pagewise]{lineno}
\usepackage{hyperref}

\begin{document}

\title[Gateaux differentiability in the space of meromorphic functions]{Gateaux differentiability in the Banach space of meromorphic functions}

\author[Mallick, Sain]{Sanjay Mallick, Debmalya Sain}
\newcommand{\acr}{\newline\indent}

\address[Mallick]{Department of Mathematics\\ Cooch Behar Panchanan Barma University\\ Cooch Behar 736101\\ West Bengal\\ INDIA}
\email{sanjay.mallick1986@gmail.com}
\address[Sain]{Department of Mathematics\\ Indian Institute of Information Technology, Raichur\\ Karnataka 584135 \\INDIA}
\email{saindebmalya@gmail.com}

\subjclass[2010]{Primary 46B20, Secondary 30D30, 30D20}
\keywords{Gateaux differentiability; Birkhoff-James orthogonality; Analytic functions; Meromorphic functions}

\begin{abstract}
	We study the Gateaux differentiability in the Banach space of meromorphic functions and obtain a complete characterization of the same, by using Birkhoff-James orthogonality techniques. We introduce the concept of extended orthogonality covering set (EOCS), which allows us to present refinements of some earlier results on the Gateaux differentiability of analytic functions. We also discuss some related properties of meromorphic functions which follow directly from the said characterization.
\end{abstract}

\maketitle

\section{Introduction}

Gateaux differentiability is undoubtedly one of the important concepts in the framework of Banach spaces, especially because of its rich connections with the analytic and the geometric properties of the underlying space. A complete description of the Gateaux differentiable points in a Banach space carries vital information on the space as well as bounded linear operators on it. We refer the readers to the pioneering study \cite{james}, and some of the later works \cite{abatzoglou, deeb-khalil, holub, RSS, sain1, sain2} for more information on this topic, from the perspective of bounded (bi)linear operators on Hilbert spaces and Banach spaces. Very recently, the Gateaux differentiability of analytic functions has been discussed in \cite{BRS}, where the authors have presented several applications of it in the study of complex analytic functions. The purpose of this article is to further refine the study conducted in \cite{BRS}, by presenting a complete characterization of the Gateaux differentiable points in the Banach space of meromorphic functions. We would like to add here that such a study may be useful in understanding the structure and properties of meromorphic functions, from a geometric perspective. Indeed, we will show that the Gateaux differentiability of a meromorphic function essentially depends on the norm attainment sets of certain functions associated naturally with the original function. This observation allows us to study a meromorphic function by using geometric techniques involving Birkhoff-James orthogonality. Although almost all the terminologies related to complex analysis, which have been used by us in this article, are standard, we refer to \cite{ahlfors} to avoid any confusion.\\

Let $f$ be a meromorphic function defined on $D:=\Gamma\cup \Int(\Gamma)$, where $\Gamma$ is a simple closed curve. Suppose that $f$ has poles at $z_{i}\in D$ of order $n_{i}$. It is clear that $f$ can have only finite number of poles in $D$, since $D$ is compact. Then we can always write $f(z)=f_{_P}(z)+f_{_R}(z)$, where $$f_{_P}(z) := \sum\limits_{z_{i}\in D}^{}\sum\limits_{j=1}^{n_{i}}\frac{a_{ij}}{(z-z_{i})^j} \;\;and\;\; f_{_R}(z) := f(z)-f_{_P}(z).$$
\par We also define  $f_{_Q}(z) := \sum\limits_{z_{i}\in D}^{}\sum\limits_{j=1}^{n_{i}}a_{ij}(z-z_{i})^j.$  
\par Let us denote by $\mathcal{M}(D)$, the vector space of all meromorphic functions defined on $D,$ under the usual operations.
\par Following \cite{cs}, define $\|.\|: \mathcal{M}(D)\rightarrow \mathbb{R}^{+}\cup\{0\}$ by 
$$\|f\| := \sup\limits_{z\in D}|f_{Q}(z)|+\sup\limits_{z\in D}|f_{R}(z)|.$$\\
 It is easy to verify that $ (\mathcal{M}(D), \|.\|) $ is a Banach space, which is the central object of our attention, as far as this article is concerned. Since the underlying norm is fixed, we will refer to $ (\mathcal{M}(D), \|.\|) $ as $\mathcal{M}(D)$ for the sake of brevity.\\
Of course, $\|f_{Q}\| = \sup\limits_{z\in D}|f_{Q}(z)|$ and $\|f_{R}\| = \sup\limits_{z\in D}|f_{R}(z)|$. In other words, $$\|f\| = \|f_{Q}\|+\|f_{R}\|.$$

As mentioned previously, our aim is to completely characterize the Gateaux differentiable points in $\mathcal{M}(D),$ with the help of Birkhoff-James orthogonality. We recall that a non-zero vector $ x $ in a normed space $ (X, \|.\|) $ is said to be Gateaux differentiable if for every vector $ y \in X, $ the following limit exists finitely:
\[ \lim_{t \to 0} \frac{\| x + ty \| - \| x \|}{t}. \]

In the language of geometry, Gateaux differentiable points are also called smooth points. It is not difficult to observe that a non-zero $ x \in X $ is a smooth point if and only if there exists a unique supporting hyperplane to the unit ball $ B_{X} $ of $ X, $ at the point $ \frac{x}{\| x \|}. $ We also recall from \cite{B}, that given any vectors $ x, y \in X, $ we say that $ x $ is Birkhoff-James orthogonal to $ y, $ written as $ x \perp_B y, $ if $ \| x + \lambda y \| \geq \| x \| $ for all scalars $ \lambda. $ Whenever the norm is induced by an inner product, the relation $ \perp_B $ is equivalent to the usual orthogonality relation $ \perp. $ However, it is also rather elementary to observe that in a general Banach space, the relation $ \perp_B $ is neither symmetric nor additive. Birkhoff-James orthogonality is extremely useful in describing many analytic properties of Banach spaces, including Gateaux differentiability. A non-zero vector $ x \in X $ is Gateaux differentiable if and only if the relation $ \perp_B $ is right-additive at $ x. $ In other words, a non-zero vector $ x \in X $ is Gateaux differentiable if and only if for every $ y, z \in X, $ the following holds:
\[ x \perp_B y,~ x \perp_B z \implies x \perp_B (y + z). \]

The above characterization will be our guiding light towards obtaining more explicit identification of Gateaux differentiable points in $ \mathcal{M}(D). $ In order to identify the Gateaux differentiable points in the Banach space of analytic functions, the following definition was introduced in \cite{BRS}.

\begin{definition}\label{d1}
A subset $ \mathcal{B} \subseteq\mathbb{C}^2$ is said to be an \textit{orthogonality covering set} (OCS, in brief), if 
\begin{align*}
\mathbb{C}= \bigcup\limits_{(z,w)\in \mathcal{B}}\{ \lambda\in \mathbb{C}: | z+\lambda w |  \geq |z| \}.
\end{align*}
\end{definition}

In order to tackle the analogous problem in the more general setting of meromorphic functions, let us now introduce the following:
\begin{definition}\label{d2}
A subset $ \mathcal{A} \subseteq\mathbb{C}^2\times\mathbb{C}^2 $ is said to be an \textit{ extended orthogonality covering set} (EOCS, in brief), if 
\begin{align*}
\mathbb{C}= \bigcup\limits_{\left((z,w),(u,v)\right)\in \mathcal{A}}\{ \lambda\in \mathbb{C}: | z+\lambda w | +  | u+\lambda v | \geq |z|+ |u| \}.
\end{align*}
\end{definition}

\par The following fact, which can be proved trivially from the definitions of EOCS and OCS, illustrates that the concept of EOCS may be regarded as an extension of the concept of OCS. 
\vspace{0.1in}\par Let  $ \mathcal{A} \subseteq\mathbb{C}^2\times\mathbb{C}^2 $ such that for any $\left((z,w),(u,v)\right)\in \mathcal{A}$ either $(z,w)=(0,0)$ or $(u,v)=(0,0)$. Define $$\mathcal{B}=\{(x,y)\in\mathbb{C}^{2}: \;\;either \left((x,y),(0,0)\right)\in \mathcal{A}\;\; or \left((0,0),(x,y)\right)\in\mathcal{A} \}.$$ Then  $ \mathcal{A}$ is an EOCS if and only if $\mathcal{B}$ is an OCS.


\vspace{0.1in}\par Given any $ f \in \mathcal{M}(D), $ we also define 
\[ M_{f_{_Q}} := \{ z \in D : | f_{_Q}(z) | = \| f_{_Q} \| \},~ M_{f_{_R}} := \{ z \in D : | f_{_R}(z) | = \| f_{_R} \| \}.  \]

Using the above concepts, first we will characterize the Birkhoff-James orthogonality in $ \mathcal{M}(D), $ in terms of EOCS. This will lead us to a tractable characterization of the Gateaux differentiable points in $ \mathcal{M}(D). $ The advantage of the obtained characterization, especially from a computational point of view, will be illustrated through several concrete examples. We will also discuss some observations on meromorphic functions, which will follow directly from our characterization.

\section{Main Results}

We begin by characterizing the Birkhoff-James orthogonality in $ \mathcal{M}(D). $

\begin{theorem}\label{t1}
Let $ f, g\in \mathcal{M}(D)$. Then the following are equivalent:\\
$ (i) $ $ f\perp_B g. $\\
$ (ii) $ $ \{ \left((f_{_Q}(z), g_{_Q}(z)), (f_{_R}(w), g_{_R}(w))\right) : (z,w)\in M_{f_{_Q}}\times M_{f_{_R}} \} $ is an EOCS.
\end{theorem}

\begin{proof}
$ (i) \implies (ii): $ Let $f\perp_B g$ and let \begin{align*}
\mathcal{A}=\{ \left((f_{_Q}(z), g_{_Q}(z)), (f_{_R}(w), g_{_R}(w))\right) : (z,w)\in M_{f_{_Q}}\times M_{f_{_R}} \}.\end{align*} Suppose on the contrary that $\mathcal{A}$ is not an EOCS. Then there exists $\lambda_{0}\in\mathbb{C}$ such that 
\begin{align*}
|f_{_Q}(z)+\lambda_{0}g_{_Q}(z)|+|f_{_R}(w)+\lambda_{0}g_{_R}(w)|<|f_{_Q}(z)|+|f_{_R}(w)|, \;\forall\; (z,w)\in M_{f_{_Q}}\times M_{f_{_R}}.\end{align*}
{\bf Claim :}  We claim that 
\begin{align*}
|f_{_Q}(z)+\mu\lambda_{0}g_{_Q}(z)|+|f_{_R}(w)+\mu\lambda_{0}g_{_R}(w)|<|f_{_Q}(z)|+|f_{_R}(w)|,\end{align*} for all $(z,w)\in M_{f_{_Q}}\times M_{f_{_R}}$ and $\mu\in(0,1)$.
\\For any fixed $(z_{0},w_{0})\in   M_{f_{_Q}}\times M_{f_{_R}}$ and  $\mu_{0}\in(0,1)$,
\beas && |f_{_Q}(z_{0})+\mu_{0}\lambda_{0}g_{_Q}(z_{0})|+|f_{_R}(w_{0})+\mu_{0}\lambda_{0}g_{_R}(w_{0})|\\&=&|(1-\mu_{0})f_{_Q}(z_{0})+\mu_{0}(f_{_Q}(z_{0})+\lambda_{0}g_{_Q}(z_{0}))|\\&&+|(1-\mu_{0})f_{_R}(w_{0})+\mu_{0}(f_{_R}(w_{0})+\lambda_{0}g_{_R}(w_{0}))|\\&\leq&(1-\mu_{0})|f_{_Q}(z_{0})|+\mu_{0}|f_{_Q}(z_{0})+\lambda_{0}g_{_Q}(z_{0})|\\&&+(1-\mu_{0})|f_{_R}(w_{0})|+\mu_{0}|f_{_R}(w_{0})+\lambda_{0}g_{_R}(w_{0})|\\&=&(1-\mu_{0})\left(|f_{_Q}(z_{0})|+|f_{_R}(w_{0})|\right)\\&&+\mu_{0}\left(|f_{_Q}(z_{0})+\lambda_{0}g_{_Q}(z_{0})|+|f_{_R}(w_{0})+\lambda_{0}g_{_R}(w_{0})|\right)\\&<& (1-\mu_{0})\left(|f_{_Q}(z_{0})|+|f_{_R}(w_{0})|\right)+\mu_{0}\left(|f_{_Q}(z_{0})|+|f_{_R}(w_{0})|\right)\\&=&|f_{_Q}(z_{0})|+|f_{_R}(w_{0})|.\eeas Hence the claim is proved.
\\Define $h: \Gamma\times\Gamma\times[0,1]\rightarrow\mathbb{R}$ by 
$$h(z,w,\mu)=|f_{_Q}(z)+\mu\lambda_{0}g_{_Q}(z)|+|f_{_R}(w)+\mu\lambda_{0}g_{_R}(w)|.$$
Then \begin{align*}
h(z_{0},w_{0},\mu_{0})<|f_{_Q}(z_{0})|+|f_{_R}(w_{0})|=\|f_{_Q}\|+\|f_{_R}\|=\|f\|.\end{align*}
Clearly, $h$ is continuous as $f_{_Q}, g_{_Q}, f_{_R}, g_{_R}$ are so. Therefore, there exist $\epsilon_{z_{0}},\epsilon_{w_{0}}$ and $\delta_{z_{0},w_{0}}$ such that \begin{align*}
h(z,w,\sigma)<\|f\|~ for~ all~ (z,w,\sigma)\in B(z_0,\epsilon_{z_{0}})\times B(w_0,\epsilon_{w_{0}})\times(\mu_{0}-\delta_{z_{0},w_{0}},\mu_{0}+\delta_{z_{0},w_{0}}).\end{align*}
In particular, \begin{align*}
h(z,w,\mu_{0})<\|f\| ~for~ all~ (z,w)\in B(z_0,\epsilon_{z_{0}})\times B(w_0,\epsilon_{w_{0}}).\end{align*}
Therefore, from the claim, we obtain that
\begin{align*}
h(z,w,\sigma)<\|f\|~ for~ all~ (z,w,\sigma)\in B(z_0,\epsilon_{z_{0}})\times B(w_0,\epsilon_{w_{0}})\times(0,\mu_{0}).\end{align*}
Suppose $(u_{0},v_{0})\in (\Gamma\times\Gamma) \setminus (M_{f_{_Q}}\times M_{f_{_R}})$. Then it is rather easy to see that
\[h(u_{0},v_{0},0) = |f_{_Q}(u_{0})| + |f_{_R}(v_{0})| < \|f_{_Q}\|+\|f_{_R}\| = \|f\|.\]
Since $h$ is continuous, there exist $\epsilon_{u_{0}},\epsilon_{v_{0}}$ and $\delta_{u_{0},v_{0}}$ such that \begin{align*}
h(u,v,\rho)<\|f\|~ for~ all~ (u,v,\rho)\in B(u_0,\epsilon_{u_{0}})\times B(v_0,\epsilon_{v_{0}})\times(-\delta_{u_{0},v_{0}},\delta_{u_{_0},v_{_0}}).\end{align*}
Therefore,
\begin{align*}
h(u,v,\delta_{u_{_0},v_{_0}})<\|f\|~ for~ all~ (u,v)\in B(u_0,\epsilon_{u_{0}})\times B(v_0,\epsilon_{v_{0}}).\end{align*}
Now proceeding similarly like the claim proved above, we can show that
\begin{align*}
h(u,v,\sigma)<\|f\|~ for~ all~ (u,v,\rho)\in B(u_0,\epsilon_{u_{0}})\times B(v_0,\epsilon_{v_{0}})\times(0,\delta_{u_{_0},v_{_0}}).\end{align*}
Since $\Gamma$ is compact, so is $\Gamma\times\Gamma$.
Now, 
\begin{align*}
\Gamma\times\Gamma\subseteq\bigcup\limits_{(z,w)\in M_{f_{_Q}}\times M_{f_{_R}}}B(z,\epsilon_{z})\times B(w,\epsilon_{w})\bigcup\bigcup\limits_{(u,v)\in \Gamma\times\Gamma\setminus M_{f_{_Q}}\times M_{f_{_R}}}B(u,\epsilon_{u})\times B(v,\epsilon_{v})
\end{align*}
By the compactness of $\Gamma\times\Gamma,$ the above open cover has a finite sub-cover, say, 
\begin{align*}
\bigcup\limits_{i=1}^{k_{1}}B(z_{i},\epsilon_{z_{i}})\times B(w_{i},\epsilon_{w_{i}})\bigcup\bigcup\limits_{j=1}^{k_2}B(u_{j},\epsilon_{u_{j}})\times B(v_{j},\epsilon_{v_{j}})
\end{align*}
where, $(z_{i},w_{i})\in M_{f_{_Q}}\times M_{f_{_R}}$ and $(u_{j},v_{j})\in \Gamma\times\Gamma\setminus M_{f_{_Q}}\times M_{f_{_R}}$.
\\\\Choose $\sigma_{0}>0$ such that $ \sigma_{0}<\mu_{z_{i},w_{i}}$ and $\sigma_{0}<\delta_{u_{j},v_{j}}$ for all $i=1,2,\ldots,k_{1}$ and for all $j=1,2,\ldots,k_2$.
Since $f_{_Q}+\sigma_{0}\lambda_{0}g_{_Q}$ and $f_{_R}+\sigma_{0}\lambda_{0}g_{_R}$ are analytic  in $D$, there exists 
\begin{align*}
(z_{1},w_{1})\in M_{f_{_Q}+\sigma_{0}\lambda_{0}g_{_Q}}\times M_{f_{_R}+\sigma_{0}\lambda_{0}g_{_R}}
~~so ~~that~~
 (z_{1},w_{1})\in \Gamma\times\Gamma.
 \end{align*}
 Now, 
 \begin{align*}
 \|f+\sigma_{0}\lambda_{0}g\|&=\sup\limits_{z\in D}|f_{_Q}(z)+\sigma_{0}\lambda_{0}g_{_Q}(z)|+\sup\limits_{w\in D}|f_{_R}(w)+\sigma_{0}\lambda_{0}g_{_R}(w)|\\&=|f_{_Q}(z_{1})+\sigma_{0}\lambda_{0}g_{_Q}(z_{1})|+|f_{_R}(w_{1})+\sigma_{0}\lambda_{0}g_{_R}(w_{1})|.
 \end{align*}
 Clearly, either $(z_{1},w_{1})\in M_{f_{_Q}}\times M_{f_{_R}}$ or $(z_{1},w_{1})\in\Gamma\times\Gamma\setminus (M_{f_{_Q}}\times M_{f_{_R}}).$ Therefore, by virtue of the choice of $\sigma_{0},$ in either case it follows that  $$|f_{_Q}(z_{1})+\sigma_{0}\lambda_{0}g_{_Q}(z_{1})|+|f_{_R}(w_{1})+\sigma_{0}\lambda_{0}g_{_R}(w_{1})|<\|f\|.$$ However, this implies that $$\|f+\sigma_{0}\lambda_{0}g\|<\|f\|,$$  a contradiction to our hypothesis  $f\perp_{B} g$. Hence $\mathcal{A}$ must be an EOCS.\\
 
 \noindent $ (ii) \implies (i): $ Assume that $\mathcal{A}$ is an EOCS. Let $\lambda\in\mathbb{C}$ be arbitrary. Then there exists $(z_{2},w_{2})\in  M_{f_{_Q}}\times M_{f_{_R}}$ such that
 $$|f_{_Q}(z_{2})+\lambda g_{_Q}(z_{2})|+|f_{_R}(w_{2})+\lambda g_{_R}(w_{2})|\geq|f_{_Q}(z_{2})|+|f_{_R}(w_{2})|.$$ 
 Now, \begin{align*}
 \|f+\lambda g\|&=\sup\limits_{z\in D}|f_{_Q}(z)+\lambda g_{_Q}(z)|+\sup\limits_{w\in D}|f_{_R}(w)+\lambda g_{_R}(w)|\\&\geq |f_{_Q}(z_{2})+\lambda g_{_Q}(z_{2})|+|f_{_R}(w_{2})+\lambda g_{_R}(w_{2})|\\&\geq |f_{_Q}(z_{2})|+|f_{_R}(w_{2})|\\&=\|f_{_Q}\|+\|f_{_R}\|\\&=\|f\|.
 \end{align*}
 Since $\lambda$ was chosen arbitrarily, it follows that $f\perp_{B} g$.
\end{proof}

In order to have a tractable characterization of the Gateaux differentiability in $ \mathcal{M}(D), $ we will depend heavily on the concept of EOCS. As it turns out, we require to identify which singleton subsets of $ \mathbb{C}^2 \times \mathbb{C}^2 $ are EOCS. We will obtain a complete description of the same, by using the Birkhoff-James orthogonality in $ l_1^2 (\mathbb{C}). $ We recall that $ l_1^2 (\mathbb{C}) $ is simply the vector space $ \mathbb{C}^2, $ equipped with the norm
\[ \| (z_1, z_2) \| = | z_1 | + | z_2 | ~\forall~ z_1, z_2 \in \mathbb{C}. \]

{\em Definition \ref{d2}} asserts that for any $ (z,u), (w,v) \in l_1^2 (\mathbb{C}), $ $ (z,u)\perp_{B}(w,v) $ if and only if \begin{align*}
	| z+\lambda w | +  | u+\lambda v | \geq |z|+ |u|~~for~~all~~\lambda\in\mathbb{C}.
\end{align*} 


\begin{rem}\label{r1}
Note that $\mathcal{A}=\{ \left((z,w),(u,v)\right)\}$ is an EOCS if and only if $(z,u)\perp_{B}(w,v)$ in $ l_1^2 (\mathbb{C}). $ Therefore, to obtain a satisfactory answer to our query which singleton subsets of $ \mathbb{C}^2 \times \mathbb{C}^2 $ are EOCS, it suffices to have an explicit description of the relation $ \perp_B $ in $ l_1^2 (\mathbb{C}). $ The following result, which completely describes the Birkhoff-James orthogonality in $ l_1^2 (\mathbb{C}), $ can be proved by following the techniques used in the proof of Theorem $ 3.3 $ of \cite{BRS2}.
\end{rem}

\begin{theorem}\label{t2}
Let $(z,u)\in l_1^2 (\mathbb{C})$ with $z\neq0, u\neq0$. Then for any $ (w,v) \in l_1^2 (\mathbb{C}), $ $(z,u)\perp_{B}(w,v)$ if and only if $$\frac{\ol z}{|z|}w+\frac{\ol u}{|u|}v=0.$$
If $z=0$ and $ u \neq 0, $ then $(z,u)\perp_{B}(w,v)$ if and only if $$aw+\frac{\ol u}{|u|}v=0,~~where~~|a|\leq1.$$
If $ z \neq 0 $ and $u=0,$ then $(z,u)\perp_{B}(w,v)$ if and only if $$\frac{\ol z}{|z|}w+bv=0,~~where~~|b|\leq1.$$
If $ z = u = 0, $ then it follows trivially that $(z,u) \perp_{B} (w,v) $ for every $(w,v)\in l_1^2(\mathbb{C}).$ 
\end{theorem}

We now present an explicit description of the singleton subsets of $ \mathbb{C}^2 \times \mathbb{C}^2 $ which are EOCS.

\begin{theorem}\label{t3}
Let $\mathcal{A}=\{ \left((z,w),(u,v)\right)\}$. Then $\mathcal{A}$ is an EOCS if and only if exactly one of the following conditions holds:
\begin{enumerate}
\item[(i)] $z=0$, $u=0$;
\item[(ii)] $z=0$, $u\neq0$ with $|w|\geq|v|$;
\item[(iii)] $z\neq0$, $u=0$ with $|w|\leq |v|$;
\item[(iv)] $z\neq0$, $u\neq0$ with $\displaystyle w=-\frac{\ol u}{|u|}\frac{|z|}{\ol z}v$;
\end{enumerate}
\end{theorem}

\begin{proof}
Let $\mathcal{A}$ be an EOCS. Suppose on the contrary that none of the conditions $ (i)-(iv) $ holds. Then one of the following conditions must hold:
\begin{enumerate}
\item[(1)] $z=0$, $u\neq0$ with $|w|<|v|$;
\item[(2)] $z\neq0$, $u=0$ with $|w|> |v|$;
\item[(3)] $z\neq0$, $u\neq0$ with $w \neq-\frac{\ol u}{|u|}\frac{|z|}{\ol z}v$;
\end{enumerate}
{\bf Case 1:} Let (1) holds. Then only the following sub-cases may occur. 
\\{\bf Case 1.1:} Let $w\neq 0$, $v\neq 0$.
Then for any $a\in\mathbb{C},$ the equation 
\begin{align*}
& aw+\frac{\ol u}{|u|}v=0 \\&\implies\frac{v}{w}=-a\frac{|u|}{\ol u}\\&\implies |a|>1,
\end{align*}
which in view of {\em Theorem \ref{t2}} implies that $(0,u)\not\perp_{B}(w,v)$ in $l_1^2(\mathbb{C})$. Hence by {\em Remark \ref{r1}}, we get that $\mathcal{A}$ is not an EOCS, a contradiction.
\\{\bf Case 1.2:} Let $w=0$, $v\neq0$. Then 
\begin{align*}
|z+\lambda w|+|u+\lambda v|= |u+\lambda v|<|0|+|u|~~for~~\lambda=-\frac{u}{2v},
\end{align*} which implies that $\mathcal{A}$ is not an EOCS, a contradiction.
\\{\bf Case 2:} Let (2) holds. Then we can have only the following sub-cases.
\\{\bf Case 2.1:} Let $w\neq0$, $v\neq0$. This case can be dealt in a manner similar to {\em Case 1.1}, to obtain a contradiction.
\\{\bf Case 2.2:} Let $w\neq0$, $v=0$. This case can be resorted similarly like {\em Case 1.2}.
\\{\bf Case 3:} Let (3) holds. Then  the following sub-cases may occur.
\\{\bf Case 3.1:}  Let $w\neq0$, $v\neq0$.  Then we have $w \neq-\frac{\ol u}{|u|}\frac{|z|}{\ol z}v$, which in view of {\em Theorem \ref{t2}} implies that $(z,u)\not\perp_{B}(w,v)$  in $l_1^2(\mathbb{C})$. Hence $\mathcal{A}$ is not an EOCS, a contradiction.
\\{\bf Case 3.2:} Let $w\neq0$, $v=0$. Then 
\begin{align*}
|z+\lambda w|+|u+\lambda v|=|z+\lambda w|+|u|<|z|+|u|~~for~~\lambda=-\frac{z}{2w},
\end{align*} which implies that $\mathcal{A}$ is not an EOCS, a contradiction.
\\{\bf Case 3.3:} Let $w=0$, $v\neq0$. This case can be resolved similarly as {\em Case 3.2}.
\\ Therefore,  for $\mathcal{A}$ to be an EOCS, one of the conditions $ (i)-(iv) $ must hold.
\vspace{0.1in}\par Conversely, let one of the conditions $ (i)-(iv) $ holds.
\\{\bf Case 1:} Let $ (i) $ holds. This case is trivial.
\\{\bf Case 2:} Let $ (ii) $ holds. Then only the following sub-cases may occur.
\\{\bf Case 2.1:} Let $w\neq0$, $v\neq0$ and $a=-\frac{v}{w}\frac{\ol u}{|u|}$. Obviously, $|a|\leq 1$ and $aw+\frac{\ol u}{|u|}v=0$. Hence, in view of {\em Theorem \ref{t2}}, we have $(0,u)\perp_{B}(w,v)$, which implies that $\mathcal{A}$ is an EOCS.
\\{\bf Case 2.2:} Let $w\neq0$, $v=0$. This case is trivial.
\\{\bf Case 2.3:} Let $w=0$, $v=0$. This case is also trivial.
\\{\bf Case 3:} Let $(iii) $ holds. Then we have the following sub-cases.
\\{\bf Case 3.1:} Let $w\neq0$, $v\neq0$ and $b=-\frac{w}{v}\frac{\ol z}{|z|}.$ Then $|b|\leq 1$ and $\frac{\ol z}{|z|}w+bv=0$. Hence, in view of {\em Theorem \ref{t2}} and {\em Remark \ref{r1}}, it follows that $ \mathcal{A}$ is an EOCS.
\\{\bf Case 3.2:}  Let $w=0$, $v\neq0$. This case is trivial.
\\{\bf Case 3.3:} Let $w=0$, $v=0$. This case is trivial.
\\{\bf Case 4:} Let $ (iv) $ holds. Then we can have the following sub-cases.
\\{\bf Case 4.1:} Let $w\neq0$, $v\neq0$. This case follows from {\em Theorem \ref{t2}} and {\em Remark \ref{r1}}.
\\{\bf Case 4.2:} Let $w=0$, $v=0$. This case is trivial.
\end{proof}

\begin{rem}
Note that all the conditions $ (i)-(iv) $ in the statement of {\em Theorem \ref{t3}} are mutually exclusive. Further observe that the conditions $ (1)-(3) $ in the beginning of the proof of this theorem, along with the conditions $ (i)-(iv) $ in the statement of {\em Theorem \ref{t3}} cover all possible combinations of $z,u,w,v.$
\end{rem}
\vspace{0.1in}\par We are now ready to present the main result of this article, which allows us to have a tractable characterization of the Gateaux differentiable points in $ \mathcal{M}(D). $

\begin{theorem}\label{t4}
Let $f\in\mathcal{M}(D)$. Then $f$ is  Gateaux differentiable in $\mathcal{M}(D)$ if and only if $M_{f_{_Q}}\times M_{f_{_R}}$ is a singleton. 
\end{theorem}

\begin{proof}
Let us first assume that $f$ is Gateaux differentiable in $\mathcal{M}(D)$. Since $f_{{Q}}$ and $f_{_R}$ are analytic in $D$, it follows that $M_{f_{_Q}}$ and $M_{f_{_R}}$ are non-empty. Suppose on the contrary that $M_{f_{_Q}}$ contains at least two elements, say, $z_{1},w_{1}.$ Let $ z_2 \in M_{f_{_R}}.$  Then \begin{align*}
(z_{1},z_{2}), (w_{1},z_{2})\in M_{f_{_Q}}\times M_{f_{_R}}, ~~~where~~ z_{1}\neq w_{1}. 
\end{align*} 
Let us now consider the following cases.
\\{\bf Case 1:} Let $f_{_Q}(z_{1})=0$. 
Suppose $f_{_R}(z_{2})=c_{1}$, for some $c_{1}\in\mathbb{C}$.
\\Define $g_{_1} : D\rightarrow\mathbb{C}$ such that \\$g_{{_1}_{Q}}(z_{1})=c_{1}$, $g_{{_1}_{R}}(z_{2})=c_{1}$ and $g_{{_1}_{Q}}(w_{1})=f_{_Q}(w_1)$.
Now, for $g_{_1}$ we have,
\begin{align*}
\{\left((f_{_Q}(z_{1}),g_{{_1}_{Q}}(z_{1})),(f_{_R}(z_{2}),g_{{_1}_{R}}(z_{2}))\right)\}=\{\left((f_{_Q}(z_{1}),c_{1}),(c_{1},c_{1})\right)\}=S_{g_{_1}} (say). 
\end{align*}
{\bf Case 1.1:} Let $c_{1}=0$. Then $S_{g_{_1}}$ falls under the condition $(i)$ of {\em Theorem \ref{t3}}. Hence in view of {\em Theorem \ref{t1}} and {\em Theorem \ref{t3}},  $f\perp_{B} g_{_1}$.
\\{\bf Case 1.2:} Let $c_{1}\not=0$. Then $S_{g_{_1}}$ falls under the condition $ (ii) $ of {\em Theorem \ref{t3}}. Hence in view of {\em Theorem \ref{t1}} and {\em Theorem \ref{t3}},  $f\perp_{B} g_{_1}$.
\\{\bf Case 2:} Let $f_{_Q}(z_{1})\not=0$. \\{\bf Case 2.1:} Suppose $f_{_R}(z_{2})=c_1=0$. Then defining $g_{_1}$ similarly like {\em Case 1}, we find that $S_{g_{_1}}$ falls under the  condition $ (iii) $ of {\em Theorem \ref{t3}}. Hence in view of {\em Theorem \ref{t1}} and {\em Theorem \ref{t3}},  $f\perp_{B} g_{_1}$.
\\{\bf Case 2.2:} Suppose $f_{_R}(z_{2})=c_1\not=0$. Then define $g_{_1} : D\rightarrow\mathbb{C}$ such that \\$g_{{_1}_{Q}}(z_{1})=-\frac{\ol { f_{_R}(z_2)}}{|f_{_R}(z_2)|}\frac{|f_{_Q}(z_{1})|}{\ol { f_{_Q}(z_{1})}}g_{_{1_{R}}(z_{2})}$, $g_{{_1}_{R}}(z_{2})=f_{_R}(z_{2})=c_{1}$ and $g_{{_1}_{Q}}(w_{1})=f_{_Q}(w_1)$. 
Hence for $g_{_1}$ we find that $S_{g_{_1}}$ falls under the condition $ (iv) $ of {\em Theorem \ref{t3}}. Therefore, in view of {\em Theorem \ref{t1}} and {\em Theorem \ref{t3}},  $f\perp_{B} g_{_1}$.
\vspace{0.1in}\par  Now we define $g_{_2} : D\rightarrow\mathbb{C}$ such that \begin{align*}
g_{{_2}_{P}}=f_{_P}-g_{{_1}_{P}}~~and~~g_{{_2}_{R}}=f_{_R}-g_{{_1}_{R}}.
\end{align*}
Therefore, $
g_{{_2}_{Q}}=f_{_Q}-g_{{_1}_{Q}}~~and~~ g_{{_2}_{Q}}(w_{1})=0,~~g_{{_2}_{R}}(z_{2})=0.
$
\\Hence we get $g_{_2}=g_{{_2}_{P}}+g_{{_2}_{R}}=f_{_P}-g_{{_1}_{P}}+f_{_R}-g_{{_1}_{R}}=f-g_{_1}$.
\\Now for $g_{_2}$, we find \begin{align*}
\{\left((f_{_Q}(w_{1}),g_{{_2}_{Q}}(w_{1})),(f_{_R}(z_{2}),g_{{_2}_{R}}(z_{2}))\right)\}=\{\left((f_{_Q}(w_{1}),0),(f_{_R}(z_{2}),0)\right)\}=S_{g_{_2}} ~(say). 
\end{align*}

Note that whatever be the values of $f_{_Q}(w_{1})$ and $f_{_R}(z_{2})$, $S_{g_{_2}}$ always falls under either one of the conditions $ (i), (ii), (iii), (iv) $ of {\em Theorem \ref{t3}}. Hence in view of {\em Theorem \ref{t1}} and {\em Theorem \ref{t3}}, $f\perp_{B} g_{_2}$. Therefore, depending upon the values of  $f_{_Q}(w_{1})$ and $f_{_R}(z_{2}),$ we can always construct $g_{_1}$ and $g_{_2}$ such that  $f\perp_{B} g_{_1}$ and $f\perp_{B} g_{_2}$ but $f\not\perp_{B} (g_{_1}+g_{_2})$.
Clearly, this contradicts our hypothesis that $ f $ is Gateaux differentiable in $ \mathcal{M}(D). $ Therefore, $M_{f_{_Q}}$ must be a singleton. Using similar techniques, it can also be shown that $M_{f_{_R}}$ is also a singleton. This completes the proof of the necessary part of the theorem.
\vspace{0.1in}\par Conversely, suppose that $M_{f_{_Q}}\times M_{f_{_R}}$ is a singleton. Observe that we may assume without any loss of generality that $M_{f_{_Q}}\times M_{f_{_R}}=\{(z_{1},z_{2})\}$ and $f_{_Q}$, $f_{_R}$ are both non-constant functions. Indeed, if one of $f_{_Q}$, $f_{_R}$ is constant, then it is trivial to see that $M_{f_{_Q}}\times M_{f_{_R}}$ is not a singleton.
\vspace{0.1in}\par   Let $f\perp_{B} g_{_1}$. Then 
$
\{\left((f_{_Q}(z_{1}),g_{{_1}_{Q}}(z_{1})),(f_{_R}(z_{2}),g_{{_1}_{R}}(z_{2}))\right)\} 
$ is an EOCS. Since $f_{_Q}(z_{1})\neq0$, $f_{_R}(z_{2})\neq 0,$ it follows from {\em Theorem \ref{t3}} that 
\begin{align*}
either~~(i)~~g_{{_1}_{Q}}(z_{1})=0~~and~~g_{{_1}_{R}}(z_{2})=0,\\or, (ii) ~~ g_{{_1}_{Q}}(z_{1})\not=0~~and~~g_{{_1}_{R}}(z_{2})\not=0\\ ~~with~~ g_{{_1}_{R}}(z_{2})=c_{2}~~ and~~g_{{_1}_{Q}}(z_{1})=-\frac{\ol { f_{_R}(z_2)}}{|f_{_R}(z_2)|}\frac{|f_{_Q}(z_{1})|}{\ol { f_{_Q}(z_{1})}}g_{_{1_{R}}(z_{2})}=c_1.
\end{align*}
\par  Let $f\perp_{B} g_{_2}$. Then as before, we get that $$
\{\left((f_{_Q}(z_{1}),g_{{_2}_{Q}}(z_{1})),(f_{_R}(z_{2}),g_{{_2}_{R}}(z_{2}))\right)\}~~ 
 is ~~an~~ EOCS~~and ~~hence,$$  
 \begin{align*}
 either~~(i)~~g_{{_2}_{Q}}(z_{1})=0~~and~~g_{{_2}_{R}}(z_{2})=0,\\or, (ii) ~~ g_{{_2}_{Q}}(z_{1})\not=0~~and~~g_{{_2}_{R}}(z_{2})\not=0\\ ~~with~~ g_{{_2}_{R}}(z_{2})=c_{4}~~ and~~g_{{_2}_{Q}}(z_{1})=-\frac{\ol { f_{_R}(z_2)}}{|f_{_R}(z_2)|}\frac{|f_{_Q}(z_{1})|}{\ol { f_{_Q}(z_{1})}}g_{_{2_{R}}(z_{2})}=c_3.
 \end{align*}
 Now, $(g_{_1}+g_{_2})_{_Q}(z_{1}) \in \{0, c_{1}, c_{3}, c_{1}+c_{3}\}$ and $(g_{_1}+g_{_2})_{_R}(z_{2}) \in \{0, c_{2}, c_{4}, c_{2}+c_{4}\}.$ Since 
 \begin{align*}
 \frac{c_{1}}{c_{2}}=\frac{c_{3}}{c_{4}}=\frac{c_{1}+c_{3}}{c_{2}+c_{4}},
 \end{align*}
it is clear that
 
 $$\{\left((f_{_Q}(z_{1}),(g_{_1}+g_{_2})_{Q}(z_{1})),(f_{_R}(z_{2}),(g_{_1}+g_{_2})_{R}(z_{2}))\right)\}$$ 
  is an EOCS.
  Therefore, $f\perp_{B} (g_{_1}+g_{_2})$. As $ g_{_1}, g_{_2} \in \mathcal{M}(D)$ were chosen arbitrarily, it follows that $f$ is Gateaux differentiable in $\mathcal{M}(D)$. This completes the proof of the sufficient part of the theorem, and establishes it completely.
\end{proof}

The following remark is particularly relevant in the context of the previous theorem:

\begin{rem}
	Although the proof of Theorem \ref{t4} depends heavily on the concept of EOCS, it should be noted that the statement of Theorem \ref{t4} is independent of EOCS. Indeed, as shown in the previous theorem, the Gateaux differentiability of an element $ f \in \mathcal{M}(D) $ can be described exclusively in terms of $  M_{f_{_Q}} $ and $ M_{f_{_R}}, $ which are inherent to the function $ f $ itself. However, to establish the above characterization of Gateaux differentiability in $ \mathcal{M}(D), $ we do require the concept of EOCS. In the opinion of the authors, it will be interesting to obtain the same characterization \emph{without using either Birkhoff-James orthogonality or the newly introduced concept of EOCS.}
\end{rem}

Through the following examples, we illustrate the computational effectiveness of Theorem \ref{t4} in identifying the Gateaux differentiable (and the non-Gateaux differentiable) functions in $ \mathcal{M}(D). $
\\\par Since a meromorphic function is either a rational function or a transcendental function, we provide examples for both of these classes of functions.

\begin{exm}Let $D:=\{z:|z|\leq r\}$ for some $r>0$ and 
$$f(z)=\frac{2z^{4}-z^{3}-8z+8}{2z^2-3z+1}.$$ This is a rational function, where $f_{Q}(z)$ and $f_{R}(z)$ attain their supremum only at $z=-r$ and $z=r$ respectively. Hence, $f(z)$ is Gateaux differentiable in $\mathcal{M}(D)$.
\end{exm}

\begin{exm}Let $D:=\{z:|z|\leq r\}$ for some $r>0$ and $$f(z)= e^{e^{24z^2+89z+711}}+\frac{720z^{2}-465z+71}{54z^{3}-51z^{2}+14z-1}.$$
 This is a transcendental meromorphic function, where $f_{Q}(z)$ and $f_{R}(z)$ attain their supremum only at $z=-r$ and $z=r$ respectively. Hence, $f(z)$ is Gateaux differentiable in $\mathcal{M}(D)$. 
\end{exm}
\begin{exm}Let $D:=\{z:|z|\leq r\}$ for some $r>0$ and $$f(z)=\frac{z^{7}+z^{5}-z^{3}+z^{2}-z+2}{z^{5}-z+1}.$$ 
 This is a rational function, where $f_{R}(z)$  attains its supremum  at $z=r,-r$. Hence, $f(z)$ is not Gateaux differentiable in $\mathcal{M}(D)$.
\end{exm}
\begin{exm}
Let $D:=\{z:|z|\leq r\}$ for some $r>0$ and $$f(z)=\displaystyle\frac{e^{\frac{\pi i}{r}z}}{e^{\frac{\pi i}{r}z}-1}.$$ Then a simple computation shows that $f_{Q}(z)$ attains its maximum at all points where $|z|=r$. Hence, this is an example of a transcendental meromorphic function which is  not Gateaux differentiable in $\mathcal{M}(D)$.  
\end{exm}

We end this article with a series of corollaries to our main result Theorem \ref{t4}, which provide additional information on the Gateaux differentiability of different sub-classes of meromorphic functions in $\mathcal{M}(D)$.

\begin{cor}\label{c1}
Let $ f(z) \in \mathcal{M}(D) $ be analytic. Then 
\begin{enumerate}
\item[(i)] $ f(z) $ is not Gateaux differentiable in $ \mathcal{M}(D).$ 
\item[(ii)] If $D:=\{z:|z-\alpha|\leq r\}$ for some $r>0$ and  $\alpha\in\mathbb{C}$, then  $\displaystyle\frac{f(z)}{z-\alpha}$ is not Gateaux differentiable in $ \mathcal{M}(D).$
\end{enumerate} 
\end{cor}
\begin{proof}
\par (i) For any analytic function $f(z)$ in $D$, clearly $f_{Q}(z)=0$ which implies  $M_{f_{_Q}}=D$. Therefore, $M_{f_{_Q}}\times M_{f_{R}}$  is not a singleton and hence the result follows from {\em Theorem \ref{t4}} .
\par (ii) Let $g(z)=\displaystyle\frac{f(z)}{z-\alpha}$. Since $f(z)$ is analytic at $\alpha$, so $f(z)$ has a Taylor series expansion about $z=\alpha,$ of the form $$f(z)=\sum\limits_{n=0}^{\infty}a_{n}(z-\alpha)^{n}.$$\par If $a_{0}\neq 0$, then $g_{_{Q}}(z)=a_{0}(z-\alpha)$. Therefore, $M_{g_{_{Q}}}=\{z:|z-\alpha|=r\}$ which is not a singleton. Hence the result.\\\par If $a_0=0$, then $g_{_{Q}}(z)=0$. Therefore, $M_{g_{_{Q}}}=D$. Hence the result.
\end{proof}
 
\begin{cor}
M\"{o}bius transformations are not Gateaux differentiable in $ \mathcal{M}(D). $
\end{cor}
\begin{proof}
 Let $f(z)=\displaystyle\frac{az+b}{cz+d}$ be a M\"{o}bius transformation in $D$. Then $f(z)$ can be represented as \bea\label{e1} f(z)&=&\frac{a}{c}+\frac{bc-ad}{c}\frac{1}{cz+d},\;\;when\;\; c\neq0;\\\nonumber&=&\frac{a}{d}z+\frac{b}{d},\;\;\;\;\;\;\;\;\;\;\;\;\;\;\;\;\;\;\;\;when\;\;c=0.\eea
From (\ref{e1}), it is obvious that in any case, one of $f_{Q}(z), f_{R}(z)$  is constant. Hence either $M_{f_{Q}}$ or $M_{f_{R}}$ becomes $D$. Therefore, $f(z)$ is not Gateaux differentiable in $\mathcal{M}(D).$
\end{proof}
\begin{cor}
Any meromorphic function $ f \in \mathcal{M}(D) $ can be written as the sum of two non-Gateaux differentiable meromorphic functions $ g, h \in \mathcal{M}(D). $ 
\end{cor}
\begin{proof}
Since analytic functions are not Gateaux differentiable in $\mathcal{M}(D)$, so is $f_{R}$. On the other hand, $f_{P}$ being a meromorphic function, we have $f_{P_{R}}=0$ and $f_{P_{P}}=f_{P}$ and hence $f_{P}$ is not Gateaux differentiable in $\mathcal{M}(D)$. Therefore, the result follows from the construction  $f=f_{R}+f_{P}$.
\end{proof}

Furthermore, we would like to present the following examples which show that the assertion $(ii)$ of {\em Corollary \ref{c1}}  may or may not hold for functions of the form $$\displaystyle\frac{f(z)}{(z-\alpha)^{n}}\;\; when\;\; n\geq2.$$
\begin{exm}
	Let $D=\{z:|z-\alpha|\leq r\}$ and $$\displaystyle f(z)=\frac{(z-\alpha+1)((z-\alpha)^{n}+1)}{(z-\alpha)^{n}},\;\; where \;\;n\geq 2.$$  
	Obviously the numerator $g(z)=(z-\alpha+1)((z˘-\alpha)^{n}+1)$ is analytic in $D$. Here $f_{R}(z)$ and $f_{Q}(z)$ both attain their supremum only at $z=\alpha+r$. Hence $f(z)$ is a Gateaux differentiable function in $\mathcal{M}(D)$.
\end{exm}
\begin{exm}
	Let $D=\{z:|z-\alpha|\leq r\}$ and $\displaystyle f(z)=\frac{\sin (z-\alpha)}{(z-\alpha)^{2}}.$ 
	\\ For this function, $f_{Q}(z)$ attains its supremum at all points on $|z-\alpha|=r$. Hence $f(z)$ is not Gateaux differentiable in $\mathcal{M}(D)$.
\end{exm}
	\section{Statements and Declarations}
 	\noindent {\bf Funding:} The research of Dr. Sanjay Mallick is supported by ``Science and Engineering Research Board, Department of Science and Technology, Government of India", under the Project File No. EEQ/2021/000316.\\
 	{\bf Conflict of interest:} Both the authors of this article declare that they do not have any conflict of interest.
 \\ {\bf Author Contributions:} Both the authors have equally contributed towards the formation of the paper.
 \\ {\bf Data availability:} Data sharing not applicable to this article as no datasets were generated or analysed during the current study.
 \section{Acknowledgement}
\par The authors would like to thank the anonymous referee for his/her careful review and nice suggestions which have increased the readability of the paper. \\\par Despite working in two different branches of Mathematics, the authors could collaborate primarily because of their friendship that started in the beautiful campus of Jadavpur University, in 2005. We are therefore elated to acknowledge the common friends from those days and their invaluable friendship that has contributed significantly to the lives and works of both of us.

\end{document}